\documentclass[11pt]{article}
\usepackage{latexsym}
\usepackage{psfig}
\newcommand{\R}{\mbox{$I\!\!R$}}

\newcommand{\gl}{{\cal L}}
\newcommand{\Pm}{{\cal PM}}
\newcommand{\fdif}[1]{C^\infty(#1)}
\newcommand{\cam}{{\cal X}}

\newcommand{\Proof}{\noindent{\em Proof: }}
\newcommand{\QED}{\hskip 10pt plus 1filll $\Box$ \vspace{10pt}}

\newtheorem{theorem}{Theorem}
\newtheorem{proposition}[theorem]{Proposition}
\newtheorem{lemma}[theorem]{Lemma}
\newtheorem{observation}[theorem]{Observation}
\newtheorem{corollary}[theorem]{Corollary}
\newtheorem{definition}[theorem]{Definition}

\begin{document}
\begin{center}{\bf\mbox{\huge The loop derivative as a curvature.}}\\
\vspace{1cm}
{\bf Martin Reiris}\\
{\it email address\tt : mreiris@cmat.edu.uy}\\
\vspace{0.3cm}
{\bf Pablo Spallanzani}\\
{\it email address\tt : pablo@cmat.edu.uy}\\
\vspace{0.3 cm}
Centro de Matem\'atica, Facultad de Ciencias, Eduardo Acevedo 1139,\\
Montevideo CP 11200, Uruguay
\end{center}
\begin{abstract}

Recently, a set of tools has been developed with the purpose of the study
of Quantum Gravity. Until now, there have been very few attempts to
put these tools into a rigorous mathematical framework. This is 
the case, for example, of the so called {\it path bundle} of a manifold.
It is well known that this topological principal bundle plays the role
of a universal bundle for the reconstruction of principal bundles
and their connections. The path bundle is canonically endowed with
a parallel transport and associated with it important types of derivatives have been considered by several authors: the Mandelstam derivative, the connection
derivative and the Loop derivative (cf. \cite{GP}).
In the present article we shall give a unified viewpoint for all
 these derivatives by developing a differentiable calculus on the 
path bundle. In particular we shall show that the loop derivative
is the curvature of a canonically defined one form that we shall called
the {\it universal connection one form}.

\end{abstract}

\section{Introduction}

Recently, a set of tools has been developed with the purpose of the study
of Quantum Gravity. Until now, there have been very few attempts to
put these tools into a rigorous mathematical framework. This is 
the case, for example, of the so called {\it path bundle} of a manifold (cf. \cite{JL1,BA,TA}).
It is well known that this topological principal bundle plays the role
of a universal bundle for the reconstruction of principal bundles
and their connections. The path bundle is canonically endowed with
a parallel transport and various important types of derivatives have been
considered by several authors: the Mandelstam derivative, the connection
derivative and the Loop derivative (cf. \cite{GP}).
In the present article we shall give a unified viewpoint for all
 these derivatives by developing a differentiable calculus on the 
path bundle which is inspired on the work of Chen \cite{Chen}.

In Section 2 we shall recall the definition of the path bundle, the loop group and the canonically defined parallel transport on the path bundle.
In section 3 we shall define the notion of differentiable space and analyze its different aspects. In particular we discus what are the differentiable functions a notion that we will use to define differential operators.
Notions of differentiable structures on the path bundle have been considered
in \cite{JL1}, but the purpose in these articles was to show that it is
possible to reconstruct a principal bundle with a connection from
the path bundle and the holonomy morphism. Our objective here is different
as it will be clear below.
Using the notion of differentiable function mentioned before we can attach
to each point in the path bundle or the loop group a ``tangent space" simply
by considering the space derivations acting on smooth functions.
This will in turn allow us to define the tangent bundle and therefore vector fields as sections of it, differential forms
and exterior derivatives.

In Section 4 we shall recall the definition of the Mandelstam derivative.
In Section 5 we shall define the {\it universal connection one form}.
This will be a differential one form in the sense of Section 3 and
we shall be able to express the Connection derivative in terms of this
universal connection.
Section 6 contains our main result. We shall prove that the loop
derivative considered in \cite{GP} is nothing but the curvature
 of the universal connection
one form. As a consequence we obtain the Bianchi identities in \cite{GP}
as the usual Bianchi identities associated with the universal connection
one form. Finally in Section 7 we shall see how to represent a particular
gauge theory using the results we obtained in the previous sections.

Finally we would like to mention that it might be possible to apply our
approach to the extended group of loops studied in \cite{TA}.

\section{Group of loops, path bundle and parallel transport}

The loop group of a manifold $M$ is defined in the following way.
Let $o$ be a fixed point in $M$, and let $L$ be the set of piecewise smooth paths $\alpha(t)$, parametrized from $[0,1]$ such that $\alpha(0)=\alpha(1)=o$.
In $L$ we consider the equivalence relation that identifies paths that differ by an orientation preserving reparametrization, then we say that two paths $\alpha$ and $\beta$ are elementary equivalent if there exists $\rho$, $\xi$, $\gamma$ such that $\alpha=\rho\cdot\xi$ and $\beta=\rho\cdot\gamma\cdot\gamma^{-1}\cdot\xi$ and we define the equivalence relation $\alpha \approx \beta$ iff there exist a sequence of paths $\alpha=\alpha_0, \alpha_1, \dots, \alpha_n=\beta$ with $\alpha_i$ elementary equivalent to $\alpha_{i+1}$.

\begin{definition}
The loop group is the quotient space of $L$ by this equivalence relation and is denoted by $\gl$.
\end{definition}

\begin{observation}
$\gl$ with the product of paths has a group structure with inverse $\alpha^{-1}$ and identity the constant path.
\end{observation}

There exist others possibilities for the definition of the equivalence relation that can lead to different loop groups~\cite{GP,AL}.

Let $PM$ denote the set of piecewise smooth paths $\alpha: [0,1] \to M$, with $\alpha(0)=o$.

\begin{definition}
The path bundle, denoted $\Pm$, is the quotient space of $PM$ by the same equivalence relation that for $\gl$.
\end{definition}

The group of loops acts on $\Pm$ by left multiplication.
As is explained in next section there exist a topology the so called Barret Topology that makes $(\Pm, M, \gl, \Pi)$ topological principal $\gl$-bundle, where $\Pi: \Pm \to M$ is the function that assigns to each path its endpoint, thus the fiber over $x$, $\Pi^{-1}(x)$ is the set of paths from $o$ to $x$ module the above mentioned equivalence relation.

There is a canonical way to define the parallel transport in this bundle, given a path $\gamma$ in $M$ with initial point $x$ and final point $y$, and an element $[\alpha]$ of the fiber of $\Pm$ over $x$, where $\alpha$ is a path going from $o$ to $x$, the parallel transport of $[\alpha]$ over $\gamma$ is $[\alpha\cdot\gamma]$, which is an element of the fiber over $y$.

\section{Differentiable structure}

In this section we define differentiable functions, tangent vectors, vector fields and differential forms for the loop group and the path bundle.

We will use the multiindex notation. A multiindex is a $n$-uple of non negative integers $\alpha=(a_1,a_2,\dots,a_n)$, $|\alpha|=a_1+\dots+a_n$ and
\[
   \partial^\alpha=\frac{\partial^{|\alpha|}}{\partial x_1^{a_1} \dots \partial x_n^{a_n}},
\]
where all the derivatives are taken in $x_i=0$.

\subsection{Differentiable spaces}

\begin{definition}
A {\em differentiable space} is a set $A$ with a family of functions from open subsets of $\R^n$ to $A$, called {\em plots}, such that if $\Phi:U \to A$ is a plot and if $g:V \to U$ is differentiable, $U \subset \R^n$, $V \subset \R^k$ then $\Phi\circ g:V \to A$ is a plot.
\end{definition}

In $A$ we consider the topology induced by the plots, that is a subset 
$U \subset A$ is open if and only if $\Phi^{-1}(U)$ is open for every plot $\phi$.

A differentiable manifold $M$ is a differentiable space considering all the differentiable functions from open subsets of $\R^n$ to $M$ as plots, and in this case all the constructions that follow (tangent spaces, forms, etc) coincide with the standard ones of differential geometry.

Now we will give $\gl$ an $\Pm$ a structure of differentiable space.

\begin{definition} Let $U$ be an open set of Euclidean space.
We define a {\em homotopy}~\cite{BA} as a function  $\Phi: U \to \Pm \hbox{ (or $\gl$)}$ such that $\Phi(x)=[t \mapsto \phi(x,t)]$ were
\[ \phi: U \times [0,1] \to M \] is continuous and there exists a partition of $[0,1]$, $0=i_0<i_1<\dots<i_n=1$ such that
\[ \phi|_{U \times [i_k, i_{k+1}]}\; \hbox{ is smooth }\; \forall k=0,1,\dots,n-1. \]
\end{definition}

Then $\gl$ and $\Pm$ are differentiable spaces taking the homotopies as plots, and the topology induced by the homotopies is the Barret topology.

Given two differentiable spaces $A$ and $B$ we can give to their product $A\times B$ a structure of differentiable space; the plots in $A\times B$ are constructed as $(\Phi,\Phi')\circ g : V \to A\times B$  where $\Phi:U \to A$, $\Phi':U' \to B$ are plots and $g: V \to U\times U'$ is any differentiable function.

\subsection{Differentiable functions}

\begin{definition}
A function $f:A \to \R$, where $A$ is a differentiable space, is said to be {\em differentiable} if $f\circ\Phi: U \to \R$ is differentiable for every plot $\Phi$. The space of differentiable functions from $A$ to $\R$ is denoted by $\fdif{A}$.
\end{definition}

\begin{definition}
A function $f: U \to A$, $U \subset \R^n$ is {\em differentiable} iff is a plot.
\end{definition}

Combining the two definitions we can define that a function $f:A \to B$ is differentiable iff $f \circ \Phi$ is a plot for every plot $\Phi$, for example the product of loops $\cdot:\gl \times \gl \to \gl$ and the action of $\gl$ over $\Pm$, $\cdot: \gl \times \Pm \to \Pm$ are differentiable.

\subsection{Tangent vectors}

\begin{definition}\label{eldifop}
Let $x_0$ be a point of $A$, an {\em elemental differential operator at $x_0$} is a linear transformation $D: \fdif{A} \to \R$ such that there exists a plot $\Phi:U \to A$, $\Phi(0)=x_0$ and a multiindex $\alpha$ such that
\[
   D(f)=\partial^\alpha (f \circ \Phi), \quad \forall f \in \fdif{A}.
\]
\end{definition}

\begin{definition}
The space of {\em differential operators at $x_0$} is the vector space generated by the elemental differential operators.
\end{definition}

Thus a differential operator at $x_0$ is a finite linear combination of elemental differential operators at $x_0$.

\begin{definition}
Let $m_{x_0}$ be the ideal in $\fdif{A}$ of the functions that vanish in $x_0$, we define the {\em order} of a differential operator at $x_0$, $D$, as the minimum $n$ such that $D|m_{x_0}^{n+1}=0$
\end{definition}

\begin{observation}
All the differential operators at $x_0$ have finite order.
\end{observation}

\Proof
It suffices to see that for an elemental diff.\ op.\ at $x_0$ the order is $\leq |\alpha|$ where $\alpha$ is the multiindex in Definition~\ref{eldifop}.
\QED

\begin{definition}
The {\em tangent space} of $A$ at the point $x_0$ is the space of first order differential operators at $x_0$.
\end{definition}

We can obtain a characterization of tangent vectors from the following lemma.

\begin{lemma}
Let $D \in \hom(\fdif{A}, R)$, the next are equivalent:
\begin{enumerate}
\item $D$ vanishes on constant functions and restricted to $m_{x_0}^2$.
\item $D(fg)=D(f)g(x_0)+f(x_0)D(g)$.
\end{enumerate}
\end{lemma}

\Proof \\
($\Rightarrow$) $D((f-f(x_0))(g-g(x_0)))=0$ because $(f-f(x_0))(g-g(x_0)) \in m_{x_0}^2$; and $D((f-f(x_0))(g-g(x_0))=D(fg)-f(x_0)D(g)-D(f)g(x_0)+D(f(x_0)g(x_0))$ and $f(x_0)g(x_0)$ is a constant function.  \\
($\Leftarrow$) If $f(x_0)=g(x_0)=0$ then $D(fg)=0$, and $D(1)=D(1\cdot 1)=D(1)+D(1)$, then $D(1)=0$, then $D$ vanishes on any constant function.
\QED

\begin{corollary}\label{opdifder}
A differential operator at $x_0$ $D$ is of first order iff $D(fg)=D(f)g(x_0)+f(x_0)D(g)$.
\end{corollary}

It can be show that the tangent vectors satisfy:
\begin{enumerate}
\item If $D \in T_xA$ and $f: V \times A \to \R$, $V \subset \R^n$ is differentiable, then $D(f): V \to \R$ is differentiable, where $D(f)(y)=D(z \mapsto f(y,z))$ is $D$ applied to $f(y, z)$ leaving $y$ fixed and considering as a function only in $z$.

\item  If $g:A \times B \to \R$ is differentiable then $D_1^{(1)} D_2^{(2)} g =D_2^{(2)} D_1^{(1)} g$ where the superscript in $D^{(n)}$ means that $D$ is applied with respect to the $n$-th variable, leaving the other fixed, that is $D^{(1)}(g)(y)=D(z \mapsto g(z,y))$.

\item If $D \in T_xA$ and $g: A \times A \to \R$ is differentiable then 
\[
	D(g\circ\Delta)=D^{(1)}(g)(x)+D^{(2)}(g)(x)
\]
where $\Delta: A \to A \times A$ is the function $\Delta(x)=(x,x)$.
\end{enumerate}

1) implies that if $g: A \times B \to \R$ is differentiable then $D^{(2)}(g):A \to \R$ is differentiable because if $\Phi$ is a plot then $D^{(2)}(g)(\Phi(y))=D(h)(y)$ were $h(y,x)=g(\Phi(y),x)$, then 2) and 3) make sense.

The proof of 1) and 2) is easy from the Definition \ref{eldifop}, to prove 3) we need the following theorem.

\begin{theorem}
$T_{({x_0},{y_0})} A \times B= T_{x_0} A \oplus T_{y_0} B$.
\end{theorem}

\Proof
If $D$ is an elemental diff.\ op.\ at $(x_0,y_0)$ we can see from the definition that $D$ can be written as
\[
D=\bar D_0^{(1)} + \hat D_0^{(2)} +
\sum_{i=1}^n \bar D_i^{(1)} \hat D_i^{(2)},
\]
where $\bar D$ are diff.\ op.\ at $x_0$ and $\hat D$ are diff.\ op.\ at $y_0$, then any differential operator can be written in this way.

Next we impose that the operator $D$ be of first order, considering the cases $g(x,y)=f(x)$ and $g(x,y)=h(y)$ we see that $\bar D_0$ and $\hat D_0$ are of first order, then $D'=\sum_{i=1}^n \bar D_i^{(1)} \hat D_i^{(2)}$ also is of first order.

We want to prove that $D'$ is null, by absurd, let us assume that it is not null.

The operator $D'$ can be written in many ways as a sum $\sum_{i=1}^n \bar D_i^{(1)} \hat D_i^{(2)}$, we choose those with less number of terms, then $\bar D_i$ are linearly independent (if not we can combine terms to obtain a sum with less number of terms).

$\hat D_1$ is not null (if not we can eliminate the term $\bar D_1^{(1)} \hat D_1^{(2)}$ of the sum) there is a function $g$ that vanishes in $y_0$ such that $\hat D_1 (g) \ne 0$; for all functions $f$ that vanish in $x_0$ we define $h(x,y)=f(x)g(y)$ then $D'(h)=0$ because $h$ is a product of two functions that vanish in $(x_0,y_0)$ and $D'$ is of first order, then
\[
\bar D_1(f) \hat D_1(g) + \cdots + \bar D_n(f) \hat D_n(g)=0.
\]
for all $f$ and $\hat D_1(g) \ne 0$, then $\bar D_i$ are linearly dependent which is absurd.
\QED

Let us see now the proof of part 3). 
If $D \in T_{x_0}A$, then the operator $D'$ defined as $D'(g)=D(g \circ \Delta)$ is in $T_{(x_0,x_0)} A \times A$, then $D'=D^{(1)}_1 + D^{(2)}_2$, considering the cases $g(x,y)=f(x)$ and $g(x,y)=h(y)$ we see that $D_1=D_2=D$, this proves~3).

Given a differentiable function $f:A \to B$ where $A$ and $B$ are differentiable spaces, then we can define its differential as a map from $T_xA$ to $T_{f(x)}B$, let $D \in T_xA$, the operator $d_xf(D)$ is defined as $d_xf(D)(g)=D(g \circ f)$, where $g$ is any function on $\fdif{B}$.

We can see that if $\Phi: U\subset \R^n \to A$ is a plot with $\Phi(0)=x$ then $d_0\Phi(v)$ is a tangent vector to $A$ at $x$.

\subsection{Tangent Bundle}

We can define the tangent bundle $TA$ of a differentiable space $A$ as the disjoint union of the tangent spaces, and it can be given a differentiable structure, for this we take as an auxiliary construction the differential bundle $DA$ defined as the disjoint union of the spaces of differential operators at $x$ for every $x \in A$, and the projection $\pi:DA \to A$ such that $\pi(D)=x$ iff $D$ is a differential operator in $x$.

Given a plot $\Phi:V \times U \to A$ we define an elemental plot $\Psi:U \to DA$ as $\Psi=\partial^\alpha \Phi$ where the derivation is taken only respect to the first variables (those who lie in $V$), that is $\Psi(x)$ is the differential operator in $\Phi(0,x)$ such that 
$\Psi(x)(f)=\partial_y^\alpha f(\Phi(y,x))$.
Then we define the plots in $DA$ as finite linear combinations of elemental plots $\Psi=a_1\Psi_1+\dots+a_n\Psi_n$ where $\Psi_i: U \to DA$ are a elemental plots and $\pi\circ\Psi_i=\pi\circ\Psi_j$ in this way $\Psi_i(x)$ are differential operators over the same point, the sum is a differential operator in that point.
$TA \subset DA$ thus we define a plot in $TA$ as a function $\Phi: U \to TA$ that is a plot in $DA$.

We now investigate when the tangent bundle is locally trivial, that is for every point $x\in A$ there is a neighborhood $U$ of $x$ such that $TU$ is diffeomorphic to $U\times T_x A$ with a diffeomorphism that is linear restricted to each tangent space,
a sufficient condition for this is the existence for every point $x \in A$ of a neighborhood $U$ of $x$ and a differentiable function $f:U\times U \times U \to A$ such that:
\begin{itemize}
\item f(x,y,y)=x.
\item f(x,x,y)=y.
\end{itemize}
because we can construct the diffeomorphism $\phi:TU \to U\times T_xA$ as 
\[
\phi(D)=d_{z=\pi(D)} f(x,\pi(D),z) (D).
\]

In the case of $\gl$ we construct the function $f$ simply as $f(\alpha,\beta,\gamma)=\alpha\cdot\beta^{-1}\cdot\gamma$. In $\Pm$ we define $f(\alpha,\beta,\gamma): \Pi^{-1}(U)\times\Pi^{-1}(U)\times\Pi^{-1}\to \Pm$ where $U$ is a convex neighborhood of the final point of $\gamma$ as 
$f(\alpha,\beta,\gamma)=\alpha\cdot\epsilon_{\alpha,\beta}\cdot\beta^{-1}\cdot\gamma\cdot\epsilon$
where $\epsilon_{\alpha,\beta}$ is the straight line joining the end points of $\alpha$ and $\beta$, and $\epsilon$ is $\epsilon_{\alpha,\beta}$ moved to the final point of $\gamma$ as we can see in figure 1.

\begin{figure}[h]
\centering
\noindent\psfig{file=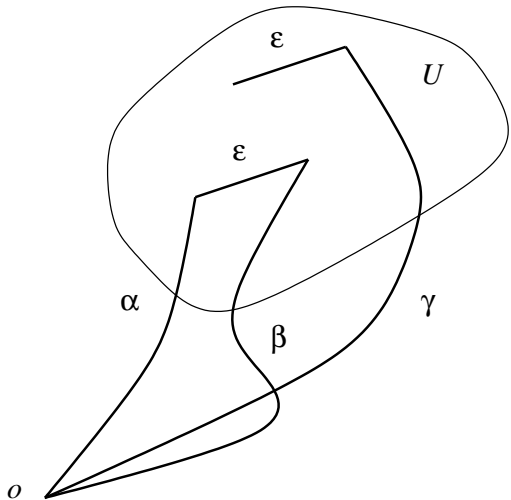}
\caption{\small Figure of $f(\alpha,\beta,\gamma)$}
\end{figure}

\subsection{Vector fields}

We define vector fields as sections of the tangent bundle, a section is a differentiable function $X: A \to TA$, such that $\pi\circ X =id$, that is $X(x)$ is a tangent vector over the point $x$. Given a section $X$ there is an associated operator $\hat X:\fdif{A} \to \fdif{A}$ such that $\hat X (f)(x)=X(x)f$.
Next we investigate this relation between sections and operators, and for technical reasons we do this in $DA$.

\begin{definition}
An operator $X: \fdif{A} \to \fdif{A}$ is an {\em elemental differential operator in $A$} iff given any plot $\Phi: U \to A$ there exists another plot $\Psi: V \times U \to A$, where $U \subset \R^k$, $V \subset \R^n$ such that
\[
   X(f) \circ \Phi = \partial^\alpha (f \circ \Psi)
\]
where the differentiation is taken respect to the firsts variables (those who lie in $V$).
\end{definition}

As an example of elemental differential operator in $\Pm$ consider the following; let $Y$ be a vector field in $M$ that vanishes in $o$, and let $\varphi_Y(x,t)$ be the flow induced by $Y$, the operator $X$ is the derivative when the path varies by the flow of $Y$, then given a homotopy $\Phi: U \to \Pm$, $\Phi(x)=[t \mapsto \phi(x,t)]$, we have the homotopy $\Psi: (-\epsilon, \epsilon) \times U \to \Pm$, $\Psi(s,x)=[t \mapsto \varphi_Y(\phi(x,t),s)]$ such that
\[
    X(f) \circ \Phi = \left. \frac{d}{ds} \right|_{s=0} f\circ\Psi.
\]

\begin{definition}
The space of {\em differential operators in $A$} is the module over $\fdif{A}$ generated by the elemental differential operators in $A$.
\end{definition}

That is, a differential operator $D$ in $A$ can be written as $D=f_1D_1+\dots+f_nD_n$, where $D_i$ are elemental differential operators and $f_i$ are differentiable functions from $A$ to $\R$.

Let $\epsilon_x$ be the operator of evaluation at $x \in A$, that is $\epsilon_x(f)=f(x)$, then if $X$ is a differential operator in $A$, $\epsilon_x \circ X$ is a differential operator at $x$.
Thus for every differential operator we can define a section $\bar X:A \to DA$ such that $\bar X(x)=\epsilon_x\circ X$, and it can be show that the operator associated with a section is always a differential operator thus we have a one to one correspondence between differential operators and sections of $DA$. In what follows we will not distinguish sections and operators.

\begin{definition}
A {\em vector field} in $A$ is a section of the tangent bundle, or equivalently, a differential operator in $A$ such that $\epsilon_x \circ X$ is of first order (that is, a tangent vector) for all $x \in A$.
\end{definition}

The space of vector fields in $A$ as a module over $\fdif{A}$ will be denoted by $\cam(A)$.

\begin{observation}\label{obsopdiff}
A differential operator in $A$ is a vector field iff $X(fg)=X(f)g+fX(g)$ for all $f,g \in \fdif{A}$.
\end{observation}

\Proof Just use Corollary \ref{opdifder}. \QED

In this way a vector field $X$ in $A$ satisfies:
\begin{enumerate}
\item If $f: V \times A \to \R$, $V \subset \R^n$ is differentiable, then $X(f): V  \times A \to \R$ is differentiable, were $X(f)(y,x)=X(z \mapsto f(y,z))(x)$ is $X$ applied to $f(y, z)$ leaving $y$ fixed and considering as function only in $z$.
\item If $g: A \times B \to \R$ is differentiable then $X_1^{(1)}X_2^{(2)}(g)=X_2^{(2)}X_1^{(1)}(g)$.
\item If $g: A \times A \to \R$ is differentiable then 
\[
	X(g\circ\Delta)=X^{(1)}(g)\circ\Delta+X^{(2)}(g)\circ\Delta
\]
\end{enumerate}

\begin{lemma}\label{opdifflemma1}
If $X$ and $Y$ are differential operators in $A$ then so is $XY$.
\end{lemma}

\Proof
It suffices to prove it when $X$ and $Y$ are elemental differential operators.

Given a homotopy $\Phi: U \to A$ then there exist $\Psi: V \times U \to A$ such that
\[
   X(Y(f)) \circ \Phi = \partial^\alpha (Y(f) \circ \Psi)
\]
because $X$ is an elemental diff.\ op., and there exists another homotopy $\Psi':V'\times V \times U \to A$ such that
\[
   Y(f) \circ \Psi = \partial^\beta (f \circ \Psi')
\]
then
\[
  XY(f) \circ \Phi = \partial^\alpha (\partial^\beta (f \circ \Psi')) =
   \partial^\gamma (f \circ \Psi')
\]
where $\partial^\alpha$ is a derivative in $V$, $\partial^\beta$ is a derivative in $V'$, and $\partial^\gamma$ is the composition of both in $V' \times V$.
\QED

\begin{definition}
The {\em bracket} is defined as the conmutator of the operators $[X_1,X_2]=X_1 X_2 - X_2 X_1$.
\end{definition}

By lemma \ref{opdifflemma1} observation \ref{obsopdiff} it can be show that the bracket of two vector fields is a vector field.

\subsection{Group of loops and Path bundle}

As mentioned before the group of loops and the path bundle become differentiable spaces taking the homotopies as plots, then the product and the inverse in $\gl$ and the action of $\gl$ over $\Pm$ are differentiable functions.

Next we define the Lie algebra of $\gl$ as the space of right invariant vector fields in $\gl$.
A vector field $X:\gl \to T\gl$ if right invariant if for every $\gamma \in \gl$ $dR_\gamma \circ X = X \circ R_\gamma$, where $R_\gamma:\gl\to\gl$ is the multiplication by $\gamma$ on the right $R_\gamma(\alpha)=\alpha\cdot\gamma$ and $dR_\gamma:T\gl \to T\gl$ is the differential of $R_g$.
The Lie bracket of two elements of the Lie algebra of $\gl$ is simply the bracket of the right invariant vector field.

\subsection{Differential forms}

The differential forms are defined as usual in differential geometry. First we define $\bigoplus^p TA$ as the Whitney sum of $p$ copies of $TA$, that is a element of $\bigoplus^p TA$ is a $p$-uple $(v_1,\dots,v_n)$ where $v_1,\dots,v_i$ are tangent vectors to $A$ over the same point.

\begin{definition}
A $p$-form $\omega$ is defined as an alternated, multilinear differentiable function $\omega:\bigoplus^p TA \to \R$.
\end{definition}

That means that for every $x$ $\omega_x:(T_xA)^p \to \R$ is a multilinear and alternated function. The differentiability of $\omega$ implies that for every $p$-uple of vector fields $X_1,\dots,X_n$ the function $x \mapsto \omega_x(X_1(x),\dots,X_p(x))$, denoted simply by $\omega(X_1,\dots,X_p)$, is differentiable.

Next we define the exterior derivative of a form, now we require that the tangent bundle of $A$ be locally trivial.
Because of local triviality of the tangent bundle every vector $v \in T_xA$ can be locally extended to a vector field $X:U \to TA$, such that $U$ is a neighborhood of $x$ and $X(x)=v$, thus it is correct to define the exterior derivative $d\omega_x(v_1,\dots,v_n)$ using vector fields $X_1,\dots,X_n$ such that $X_i(x)=v_i$ if we prove that the result only depends in $v_1,\dots,v_n$

\begin{definition}
The {\em exterior derivative} of $\omega$ is defined as:
\begin{eqnarray*}
d\omega(X_{1},\dots,X_{n+1}) &=& \frac{1}{n+1}\Bigg[\sum_{i=1}^{n+1}(-1)^{i+1}X_{i}(\omega(X_{1},\dots, \hat{X}_{i},\dots, X_{n+1}))+ \\
&& +\sum_{i<j}(-1)^{i+j}\omega([X_{i},X_{j}],X_{1},\dots, \hat{X}_{i},\dots,\hat{X}_{j},\dots, X_{n+1})\Bigg].
\end{eqnarray*}
\end{definition}

By the formula we see that $d\omega$ is also multilinear and alternated, but it remain to be proved that it is local too. Let us prove this in the next lemma.

\begin{lemma}
$d\omega(X_1,\dots,X_n)(x)$ only depends on the values $X_i(x)$.
\end{lemma}

\Proof 
Because of linearity, it is sufficient to show that $d\omega(X_{1},\dots,X_{n+1})(x)=0$ when $X_{1}(x)=0$.
Looking at the formula it is also sufficient to show that 
\[
X_{i}(x)(\omega(X_{1},\dots, \hat{X}_{i},\dots,X_{n+1}))+ \omega([X_{1},X_{i}],X_{2},\dots, \hat{X}_{i},\dots, X_{n+1})(x)=0.
\]

Let us consider the $1$-form $\omega(X)=\omega(X,X_{2},\dots,\hat{X_{i}},\dots,X_{n+1})$ then what we have to prove is 
\[
Y(x)\omega(X)+\omega([X,Y])(x)=0
\]
whenever $X(x)=0$.

Next we construct using a local trivialization of the tangent bundle a function $X:U \times U \to TA$, where $U$ is a neighborhood of $x$ that verifies
\begin{itemize}
\item $X(y,z)\in T_{y} A$
\item $X(y,y)=X(y)$
\item $X(y,x)=0$
\end{itemize}

The proof then follows looking at the equalities
\begin{eqnarray*}
Y_{y}(x)\omega(X(y,y)) &=& Y_{y}(x)\omega(X(y,x)) + Y_{y}(x)\omega(X(x,y))= \\
&=& Y_{y}(x)\omega(X(x,y))=\omega(Y_{y}(x)X(x,y)) \\{}
[X,Y](x) &=& -Y_{y}(x)(X(y))=-Y_{y}(x)(X(y,y)) =\\ 
&=& -Y_{y}(x)(X(y,x))-Y_{z}(x)(X(x,z))=
-Y_{z}(x)(X(x,z))
\end{eqnarray*}
\QED

\section{Mandelstam derivative}

The parallel transport defined in Section 3 induces a set of planes in $\Pm$ that will be the horizontal planes of the universal connection one form that we shall define later.

\begin{definition}
The Mandelstam derivative on $\pi_o^x$ in the direction of $v \in T_{x}M$ is given by
\[
D_{v} f(\pi_{o}^{x})=d_{x}(f\circ\Phi)(v)
\]
where $\phi$ is a family of curves represented by
\[
\phi(y,t)=\left\{ \begin{array}{ll}

           \pi_{o}^{x}(2t)&0\leq t\leq 1/2 \\
           g((2t-1)g^{-1}(y))  &1/2\leq t\leq 1
                          \end{array} \right.
\]
and 
\[
\Phi(y)=[t \mapsto \phi(y,t)]
\]
where $g:U \to M$ is a chart with $U$ a convex neighborhood of $0$ and $g(0)=x$.
\end{definition}

\begin{observation}
If we call $\delta v$ the segment (in the chart $g$) from $x$ to $x+\epsilon v$ then
\[
  f(\pi_o^x \cdot \delta v) = f(\pi_o^x) + \epsilon D_v f + o(\epsilon)
\]
\end{observation}

Thus the Mandelstam derivative coincides with the one defined in~\cite{GP}

It needs to be proved that the definition does not depend of the chart $g$, for that we first state a lemma proved in~\cite{BA}.

\begin{lemma}\label{b1}
If $\phi: (-\epsilon, \epsilon) \times [0,1]\to M$ is a homotopy of loops where $\phi(0,t)=o$ is the constant loop and $\Phi(s)=[t \mapsto \phi(s,t)]$ then for any differentiable function  $f: \gl \to \R$,
\[
 \left. \frac{d}{ds} \right|_{s=0} f \circ \Phi = 0.
\]
\end{lemma}

\begin{proposition}
Let $\alpha: (-\epsilon, \epsilon) \to M$ be any curve with $\alpha(0)=x$ and $\dot\alpha(0)=v$, define the homotopy of paths $\Phi(s)=\pi_o^x \cdot [t \mapsto \alpha(st)]$, and let the operator $\tilde D_v$ be such that
\[
  \tilde D_v(f) = \left. \frac{d}{ds} \right|_{s=0} f \circ \Phi
\]
then $\tilde D_v = D_v$.
\end{proposition}

\Proof
Let $\alpha(s)=(\alpha_1(s),\dots,\alpha_n(s))$ in the chart $g$ and assume that in that chart $v=e_1$ where $e_1,\dots,e_n$ is the canonical basis of $\R^n$, then $\alpha_1'=1$ and $\alpha_i'=0$, $i=2,\dots,n$.

We construct the homotopy of loops $\Phi$ such that $\Phi(s)$ is the loop based at $x$ formed by the composition of the curve $\alpha$ from $x$ to $\alpha(s)$ then the straight segment (in the chart $g$) from $\alpha(s)$ to $g((\alpha_1(s),0,\dots,0))$ and then the segment to $x$.

By lemma \ref{b1} $\left. \frac{d}{ds} \right|_{s=0} f(\pi_o^x \cdot \Phi(s)) = 0$ and expanding $f(\pi_o^x \cdot \Phi(s))$ until first order in $s$ we get
\[
f(\pi_o^x \cdot \Phi(s)) = f(\pi_o^x) + s\tilde D_v(f) + \sum_{i=2}^n \alpha_i(s) D_{e_i}(f) - \alpha_1(s)D_{e_1}(f).
\]
taking limit when $s \to 0$ we get $\tilde D_v(f)=D_v(f)$.
\QED

\section{Universal connection one form}

Following which is the habitual definition of connection one form of the theory of principal bundles, we define a connection in $\Pm$, that is a $1$-form evaluated in the Lie algebra $=T_e\gl$ of $\gl$, such that is equivariant under the action of the group and is the identity over the vertical subspaces~\cite{DU}.

Given $D \in T_\pi \Pm$ and $f \in \fdif{\gl}$ we first transform $f$ into a function on  $\Pm$ given by: $g(\gamma)=f(\gamma \cdot v(\gamma) \cdot \pi^{-1})$, where $v(\gamma)$ is the segment that joins $\Pi(\gamma)$ (the final point of $\gamma$) with
$\Pi(\pi)$ (the final point of $\pi$) in a chart $h: U \to M$ with $U$ a convex neighborhood of $0$ and $h(0)=\Pi(\pi)$. We define $\delta_\pi (D)(f) = D(g)$.

\begin{lemma}
$\delta$ is a $1$-form in $\Pm$ with values in $T_e \gl$.
\end{lemma}

\Proof
First we must show that the definition does not depend of the chart $h$ chosen, for that we investigate how a vector $D \in T_\pi\Pm$ can be decomposed into a horizontal and a vertical part.
\begin{eqnarray*}
D_\gamma f(\gamma) &=& D_\gamma f(\gamma \cdot v(\gamma) \cdot \pi^{-1} \cdot \pi \cdot v(\gamma)^{-1}) \\
 &=& D_\gamma f(\gamma \cdot v(\gamma) \cdot \pi^{-1} \cdot \pi) + D_\gamma f(\pi \cdot v(\gamma)^{-1}) 
 = dU(\delta_\pi(D))(f) + D_v(f)
\end{eqnarray*}
where $v=d\Pi(D)$ and $U: \gl \to \Pm$, $U(\xi)=\xi \cdot \pi$;
because $D_v$ does not depend of the chart we see that $\delta_\pi$ does not depend of the chart.

To verify that $\delta$ is a $1$-form in $\Pm$ evaluated in $T_e \gl$ we define the function $g: \Pm \times \Pm \to \gl$, $g(\gamma, \pi)=f(\gamma \cdot v \cdot \pi^{-1})$, then $\delta_\pi(D)(f)=D^{(1)} g(\pi, \pi)$. The same expression is written as $\delta(X)(f)=X^{(1)}(g) \circ \Delta$ and with it we check that it is differentiable.
\QED

Let us see that $\delta$ is really a connection. 

\begin{lemma}
$\delta$ satisfies the definition of a connection.
\end{lemma}

\Proof
To prove this we have to see first that it is the identity over the vertical vectors and must verify the compatibility condition over the action of the group.

Let $U$ be $U : \gl \to \Pm$, $U(\gamma)=\gamma \cdot \pi$, the vertical vectors in $T_\pi PM$ are those that lie at the image of $dU$, and making the calculation $\delta_\pi \circ dU$ we obtain
\[
\delta_\pi(dU(D))(f)=dU(D)_\gamma f(\gamma \cdot v \cdot \pi^{-1})=
D_{\gamma'} f(\gamma' \cdot \pi \cdot v \cdot \pi^{-1})=D_{\gamma'} f(\gamma')=D(f).
\]

Let $U$ be $U_\alpha : \Pm\to \Pm$, $U_\alpha (\pi)=\alpha \cdot \pi$,
\[
(U_\alpha^\ast\delta)_\pi(D)(f)=\delta_{\alpha \cdot \pi} (dU_\alpha(D))(f) = dU_\alpha(D)_\gamma f(\gamma \cdot v \cdot \pi^{-1} \cdot \alpha^{-1})
\]\[
= D_{\pi'} f(\alpha \cdot \pi' \cdot v \cdot \pi^{-1} \cdot \alpha^{-1}) = Ad(\alpha) \delta_\pi(D) (f).
\]
\QED

As an illustrative example let us see how the connection derivative given in \cite{GP} enters into this context. Let us take a section of the bundle around the point of the fiber over $x$, $\pi_{o}^{x}$. That means an election of a family of curves $\Phi : U\to \Pm$ where $U$ is an open neighborhood of the $x$ point y $\Pi(\Phi(y))=y,\ y\in U$. 
Given a function $f(\pi_{o}^{x})$ we define the tangent vector to the section $\Phi$ as the operator that applied to $f$ gives 
\[
\tilde{D}_{v}f(\pi_{0}^{x})=\partial_{v}f\circ\Phi(x)=\frac{df(\pi_{0}^{x+\epsilon v})}{d\epsilon}.
\]
We observe that the projection of this vector over $M$ is $v$.

\begin{figure}[h]
\centering
\noindent\psfig{file=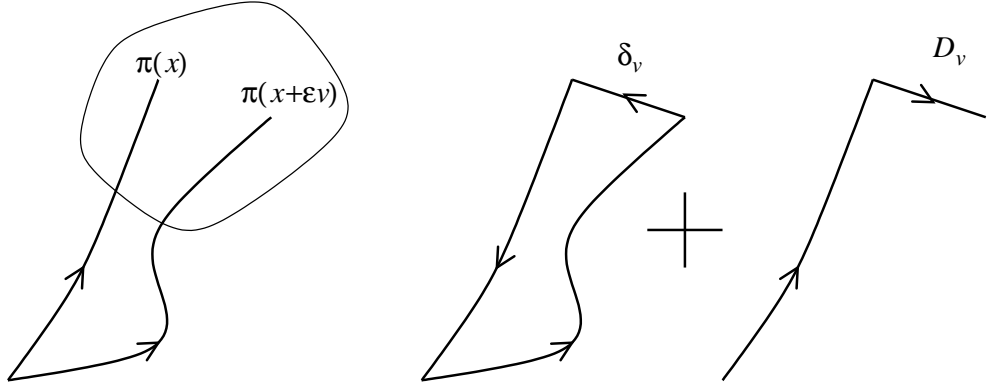}
\caption{\small The figure on the left represent a trivialization. Those on the right are the vertical and horizontal vectors decomposing the tangent vector to the trivialization.}
\end{figure}

Let us see how is the decomposition of $\tilde{D}_{v}$ in a vertical and a horizontal vector.
We write $\tilde{D}_{v}$ as
\[
\tilde{D}_{v}f(\pi_{0}^{x})=\frac{df(\pi_{0}^{x+\epsilon v}.-\epsilon v. \pi_{x}^{0}).(\pi_{0}^{x}.\epsilon v)}{d\epsilon}=v^{\mu}D_{\mu}f(\pi_{0}^{x})+
v^{\mu}\delta_{\mu}(x)f(\pi_{0}^{x})
\]
where $\delta_{\mu}$ is
\[
\delta_{\mu}f(\pi_{0}^{x})=\frac{df(\pi_{o}^{x+\epsilon v}.-\epsilon v)}{d\epsilon}.
\]

This is exactly the connection derivative given in \cite{GP} and in our context is written as
\[
\delta_{\mu}f(\pi_{o}^{x})=\delta_{\pi_{o}^{x}}(\tilde{D}_{v})_{\alpha}f(\alpha\cdot\pi_{o}^{x}).
\]

\section{The loop derivative as a curvature}

In this section, we will see that the loop derivative defined in \cite{GP} represents the curvature of the universal connection one form defined above.

\begin{definition}
The loop derivative is a vector field in $\gl$ that applied to a function $f(\gamma)$, $\gamma \in \gl$ is
\[
\Delta_{u,v} (\pi_{o}^{x}) f(\gamma)=\frac{\partial^{2} f(\pi\cdot\Box\cdot\pi^{-1}\cdot\gamma)}{\partial \epsilon_{1}\partial \epsilon_{2}}
\]
where $\Box=\Box_{\epsilon_{1},\epsilon_{2}}$ is the parallelogram, taken in a local chart of $M$, with vertex at $x$ and edges in the directions of the vectors $u$ and $v$ with lengths $\epsilon_{1}\cdot \| u\|$ and $\epsilon_{2}\cdot \| v\|$ respectively.
\end{definition}

The curvature form is written in terms of the connection form as (cf. \cite{DU})
\[
\Omega=d\delta+\frac{1}{2}[\delta,\delta]
\]

Since the curvature form is horizontal it suffices to evaluate it on horizontal vectors, that is Mandelstam derivatives $D_\mu$, $D_\nu$ that are horizontal vector fields defined in a neighborhood of the point $x$.

In this context applying the expression for the exterior derivative given in Section 3, and evaluating the curvature on the horizontal fields $D_\mu$, $D_\nu$ we obtain
\begin{eqnarray*}
\Omega(D_{\mu},D_{\nu})&=&\frac12 \Big[ D_{\mu}\delta(D_{\nu}) - D_{\nu}\delta(D_{\mu}) - \delta([D_\mu, D_\nu]) + [\delta(D_{\mu}), \delta(D_{\nu})] \Big]= \\ &=&-\frac12\delta([D_{\mu},D_{\nu}])=-\frac12[D_{\mu},D_{\nu}],
\end{eqnarray*}
because the connection form is null on horizontal vectors and is the identity on vertical vectors and $[D_\mu, D_\nu]$ is vertical.

In \cite{GP} it is proved that
\[
[D_{\mu},D_{\nu}]=\Delta_{\mu,\nu},
\]
thus
\[
\Delta_{\mu,\nu}=-2\Omega(D_\mu,D_\nu).
\]
This also proves that the loop derivative does not depend of the chart taken in the definition.

In this context the Bianchi identity presented in \cite{GP} follows naturally from the Bianchi identity for the curvature in $\Pm$
\[
d\Omega=[\delta,\Omega],
\]
\[
d\Omega(D_{\mu},D_{\nu},D_{\xi})=D_{\mu}\Delta_{\nu,\xi} + D_{\nu}\Delta_{\xi,\mu}+ D_{\xi}\Delta_{\mu,\nu}=0
\]
because $[\delta,\Omega](D_{\mu},D_{\nu},D_{\xi})=0$.

\section{Application to the representation of a particular gauge theory}

Suppose now that we want to represent a particular gauge theory in this context, that is taking the pull back of the theory to the path bundle.

It is well known that the path bundle has the universal property that every principal $(GM,M,\Pi,G)$ bundle can be reconstructed from the holonomy morphism $H:\gl \to G$ coming from the connection $\theta$ in $GM$ (cf. \cite{BA}).
 In this way every principal bundle can be seen as an extension of the path bundle~\cite{DU}.
In particular the bundle morphism is expressed

\[f(\pi_{o}^{x})=\Pi(\pi G_{o}^{x}).\]

Where $\Pi(\pi G_{o}^{x})$ is the final point of the horizontal path in $GM$, that is projected over $\Pi_{o}^{x}$. We also have the following identity between the connections

\[f^{*} \theta = d H \circ \delta.\]

Helped with this formula, the identities presented in \cite{GP} 

\[\delta_{\mu}(x)H(\gamma )=A_{\mu}(x)H(\gamma);\]

\[\Delta_{\mu,\nu}(\pi_{o}^{x})H(\gamma)=F_{\mu\nu}(x)H(\gamma),\]

are expressed as 

\[\delta_{\mu}(x)H(\gamma )=d H \circ \delta(\tilde{D}_{\mu});\]

\[\Delta_{\mu,\nu}(\pi_{o}^{x})H(\gamma)= d H\circ \Omega (D_{\mu},D_{\nu}).\]

\section{Aknowledgement}

We want to thanks specially to proffesors Gabriel Paternain, Miguel Paternain, and Rodolfo Gambini for their comments and important suggestions.

\end{document}